\documentclass{elsart}
\usepackage{amscd,amssymb,graphics}
\usepackage{mathrsfs} 
\newtheorem{theorem}{Theorem}[section]
\newtheorem{corollary}[theorem]{Corollary} 
\newtheorem{lemma}[theorem]{Lemma}
\newtheorem{proposition}[theorem]{Proposition}
\theoremstyle{definition}

\newtheorem{definition}[theorem]{Definition}
\theoremstyle{remark}

\newtheorem{example}[theorem]{Example}

\newcommand{\abs}[1]{\vert#1\vert}

\def\R {{\mathbb R}}

\def\I {{\mathbb I}}
\def\N{{\mathbb N}}
\def\Ur{{\mathbb U}}

\def\e{{\varepsilon}}

\def\Z {{\mathbb Z}}
\def\s{{\mathbb S}}

\def\Iso{{\mathrm{Iso}\,}}
\def\ev{{\mathrm{ev}\,}}

\def\Homeo{{\mathrm{Homeo}\,}}
\def\diam{{\mathrm{diam}\,}}

\begin{document}
\begin{frontmatter}

\title{The isometry group of the Urysohn
space as a L\'evy group}

\author{Vladimir Pestov}
\address{Department of Mathematics and Statistics, 
University of Ottawa, 585 King Edward Ave., Ottawa, Ontario, Canada K1N 6N5
}
\ead{vpest283@uottawa.ca}
\ead[url]{http://aix1/uottawa.ca/$^\sim$vpest283}

\begin{abstract} We prove that the isometry
group $\Iso(\Ur)$ of the universal Urysohn metric space $\Ur$ equipped with
the natural Polish topology is a L\'evy
group in the sense of Gromov and Milman, that is, admits an approximating
chain of compact (in fact, finite) subgroups, exhibiting the phenomenon of
concentration of measure. This strengthens an earlier result by Vershik stating
that $\Iso(\Ur)$ has a dense locally finite subgroup. 
\end{abstract}
\end{frontmatter}


\section{Introduction}

The following concept, introduced by P.S. Urysohn \cite{U25,U27}, has
generated a considerable and steadily growing interest over the past 
two to three decades.

\begin{definition} 
\label{def:urysohn}
The \index{Urysohn metric space} \index{U@$\Ur$}
{\em Urysohn metric space} $\Ur$ is defined by three
conditions:
\begin{enumerate}
\item $\Ur$ is a complete separable metric space;
\item $\Ur$ is ultrahomogeneous, that is, every isometry between
two finite metric subspaces of $\Ur$ extends to a global isometry of
$\Ur$ onto itself;
\item $\Ur$ is universal, that is, contains an isometric copy of every
separable metric space.
\end{enumerate}
\end{definition}

An equivalent property distinguishing $\Ur$ among complete separable metric
spaces is \textit{finite injectivity:} if $X$ is a metric subspace of a
finite metric space $Y$, then every metric embedding of $X$ into $\Ur$ extends
to a metric embedding $Y\hookrightarrow\Ur$. Establishing an equivalence
between this description and Definition \ref{def:urysohn} is an
enjoyable exercise.

Such a metric space $\Ur$ exists and is unique up to an
isometry, and in addition to the original proof by Urysohn,
there are presently several known alternative proofs of this result, most 
notably those in \cite{Kat} and in \cite{Ver98,Ver02,Ver04}.

At the same time, there is still no known
concrete realization (model) of the Urysohn space,
and finding such a model is one of the most interesting open problems of
the theory, mentioned by such mathematicians as
Fr\'echet \cite{Fre}, p. 100 and P.S. Alexandroff \cite{U72}, and presently
being advertised by Vershik. The only bit of constructive knowledge about
the structure of the Urysohn space currently available
is that $\Ur$ is homeomorphic to the
Hilbert space $\ell^2$ (Uspenskij \cite{Usp04}).

A ``poor man's version'' 
of the Urysohn space $\Ur$, the so-called \textit{random graph}
$R$ (discovered much
later than the Urysohn space, see e.g. \cite{Rado}), has a model 
(in fact, more than one, cf. \cite{Cam}). 
The random graph can be viewed as a version of the universal Urysohn metric
space whose metric only takes values $0,1,2$, which fact 
offers some hope that a model for $\Ur$ can also be found.

One can in fact study a variety of versions of the Urysohn space defined by
a restriction on the collection of values that the distance is allowed to
take. Among them, one of the most interesting and useful objects is the
\textit{Urysohn space of diameter one,} $\Ur_1$. Here the diameter of a metric
space $X$ is defined as $\sup\{d(x,y)\colon x,y\in X\}$, the space $\Ur_1$ has 
diameter $1$, and possesses the same properties as the space $\Ur$ with the
exception that it is universal for all separable metric spaces having
diameter one. It is easy to show that the space $\Ur_1$ is
isometric to the sphere of radius $1/2$ around any point in $\Ur$.

An interesting approach to the Urysohn space
was proposed by Vershik who regards $\Ur$
as a generic, or random, metric space. Here is one of his results.
Denote by $M$ the Polish space of all metrics
on a countably infinite set $\omega$. Let $P(M)$ denote the Polish
space of all
probability measures on $M$. Then, for a generic measure $\mu\in P(M)$
(in the sense of Baire category), the completion of the metric space $(X,d)$
is isometric to the Urysohn space $\Ur$ $\mu$-almost surely in $d\in M$.
We refer the reader to a very interesting theory developed 
in \cite{Ver98,Ver02} and especially \cite{Ver04}. Cf. also \cite{Ver05}.

The group of all isometries of the Urysohn space $\Ur$ onto itself,
equipped with the topology of simple convergence (or the compact-open topology, 
which happens to be the same), is a Polish (separable completely metrizable)
topological group. It possesses the following remarkable property, discovered
by Uspenskij.

\begin{theorem}
[Uspenskij \cite{Usp90}; Cf. also \cite{Gr}, 3.11.$\frac 23_+$]
\label{th:iso(u)}
The Polish group $\Iso(\Ur)$ is a universal second-countable
topological group. In other words, every second-countable topological
group $G$ embeds into $\Iso(\Ur)$ as a topological subgroup.
\qed
\label{th:univ}
\end{theorem}

The same property is shared by the topological group $\Iso(\Ur_1)$. 
Other known results include the following.

\begin{theorem}[Uspenskij \cite{Usp98}]
The group $\Iso(\Ur_1)$ is topologically simple 
(contains no non-trivial closed normal subgroups) and minimal
(admits no strictly coarser Hausdorff group topology)
\label{th:min}
\qed
\end{theorem}

One can deduce from the above fact some straightforward
but still interesting corollaries which, to this
author's knowledge, have never been stated by anyone explicitely. For
example, $\Iso(\Ur_1)$ admits no non-trivial (different from the identity)
continuous unitary representations. In fact, a stronger result holds.

\begin{corollary}
The topological group $\Iso(\Ur_1)$ admits no non-trivial continuous
representations by isometries in reflexive Banach spaces. 
\label{c:1}
\end{corollary}

\begin{pf} According to Megrelishvili \cite{Meg},
the group $\Homeo_+[0,1]$ consisting of
all orientation-preserving self-homeomorphisms of the closed unit
interval and equipped with the compact-open topology, admits no
non-trivial continuous
representations by isometries in reflexive Banach spaces. 
By Uspenskij's theorem, 
$\Homeo_+[0,1]$ embeds into $\Iso(\Ur_1)$ as a topological subgroup.
If now $\pi$ is a continuous representation of $\Iso(\Ur_1)$ in a reflexive 
Banach space $E$ by isometries, that is, a continuous
homomorphism $\pi\colon \Iso(\Ur_1)\to \Iso(E)$ where the latter group is
equipped with the strong operator topology, then, by force of
Theorem \ref{th:min}, the kernel $\ker\pi$
is either $\{e\}$ or all of $\Iso(\Ur_1)$. In the former case, the 
restriction of $\pi$ to a copy of $\Homeo_+[0,1]$ must be a continuous
faithful representation by isometries in a reflexive Banach space,
which is ruled out by Megrelishvili's theorem. We conclude that
$\ker\pi=\Iso(\Ur_1)$, that is, the representation $\pi$ is trivial (assigns the
identity operator to every element of the group).
\qed\end{pf}

Modulo a result independently obtained by
Megrelishvili \cite{Meg2} and Shtern \cite{Shtern}, this implies:

\begin{corollary}
Every continuous weakly almost periodic function on $\Iso(\Ur_1)$
is constant. \qed
\label{c:2}
\end{corollary}

An action of a topological group $G$ on a finite measure space
$(X,\mu)$ is called \textit{measurable,} or a {\em near-action}, if 
for every $g\in G$ the motion $X\ni x\mapsto gx\in X$
is a bi-measurable map defined $\mu$-almost everywhere, 
and for every measurable set
$A\subseteq X$ the function $G\ni g\mapsto \mu(gA\Delta A)\in\R$ is continuous.
In addition, the identities $g(hx) = (gh)x$ and $ex=x$
hold for $\mu$-a.e. $x\in X$
and every $g,h\in X$. Such an action is \textit{measure class preserving} if 
for every measurable subset $A\subseteq X$ and every $g\in G$, the 
set $g\cdot A$, defined up to a $\mu$-null set, has 
measure $\mu(g\cdot A)>0$ if and only if $\mu(A)>0$. Finally, we say that
an action as above is {\em trivial} if the set of $G$-fixed points has full
measure.

\begin{corollary} The topological group
$\Iso(\Ur_1)$ admits no non-trivial measurable action on a measure
space, preserving the measure class.
\label{c:3}
\end{corollary}

\begin{pf} Indeed, every such action leads to a non-trivial strongly
continuous unitary representation via the standard construction of the quasi-regular
representation in the space $L^2(X,\mu)$, given by the formula
\[^gf(x) = \left(\frac{d(\mu\circ g^{-1})}{d\mu}\right)^{\frac 12}
f(g^{-1}x),\]
where $d/d\mu$ is the Radon-Nykodim derivative.
\qed\end{pf}

We do not know if the analogues of Theorem \ref{th:min} and Corollaries
\ref{c:1}, \ref{c:2}, \ref{c:3} hold for the Polish group $\Iso(\Ur)$.

Another example of a universal Polish group was also previously discovered by
Uspenskij \cite{Usp86}: the group $\Homeo(Q)$ of
self-homeomorphisms of the Hilbert
cube $Q=\I^{\aleph_0}$ equipped with the compact-open topology. Apparently,
$\Homeo(Q)$ and $\Iso(\Ur)$ remain to date, essentially, the only known examples 
of universal Polish groups (if one discounts 
the modifications of the latter such
as $\Iso(\Ur_1)$). As pointed out in \cite{P99},
these two topological
groups are not isomorphic between themselves. Indeed, the Hilbert
cube is topologically homogeneous, that is, the
action of
$\Homeo(Q)$ on the compact space $Q$ is transitive and therefore fixed 
point-free, cf. e.g. \cite{vM}.
At the same time, the dynamic behaviour of the
groups such as $\Iso(\Ur)$ is markedly different.

\begin{definition}
One says that a topological group $G$ is 
\index{group!extremely amenable}
{\em extremely amenable,} or has the 
\index{fixed point on compacta property}
{\em fixed point on compacta property,} if
every continuous action of $G$ on a compact space $X$ admits a fixed
point: for some $\xi\in X$ and all $g\in G$, one has $g\xi=\xi$.
\end{definition}

As first noted by Granirer and Lau \cite{GraL},
no locally compact group different from the trivial group
$\{e\}$ is extremely amenable. In fact,
until an example was constructed by Herer and Christensen in \cite{HC}, the very
existence of extremely amenable topological groups remained in doubt.
However, since Gromov and Milman \cite{GrM} proved that the unitary group
$U(\ell^2)$ of a separable Hilbert space equipped with the strong operator
topology is extremely amenable, it gradually became 
clear that the property is rather
common among the concrete ``infinite-dimensional'' topological groups.
We refer the reader to two recent articles \cite{KPT} and \cite{GP2} which
together cover most of examples of extremely amenable groups known to date.

The present author had shown in \cite{P02} that the topological group
$\Iso(\Ur)$ is extremely amenable. Consequently, it is non-isomorphic,
as a topological group, to $\Homeo(Q)$. The same is true of $\Iso(\Ur_1)$.

Vershik has demonstrated in \cite{Ver05b} that the group $\Iso(\Ur)$
contains a locally finite everywhere dense subgroup. We will give an
alternative proof of this result below in 
Section \ref{s:lfs}. This proof is a step towards
theorem \ref{th:isoulevy} which is the main result of our article. 
Before stating this result, we need to remind some concepts introduced by
Gromov and Milman \cite{GrM} and linking topological dynamics of
``large'' groups with asymptotic geometric analysis \cite{M00}. 

The 
{\em phenomenon of concentration of
measure on high-dimensional structures} says, intuitively speaking, that
the geometric structures -- such as the Euclidean spheres --
of high finite dimension typically 
have the property that an overwhelming proportion
of points are very close to every set containing at least half of the
points. Technically, the phenomenon is dealt with in the following
framework. 

\begin{definition}[Gromov and Milman \cite{GrM}]
A {\em space with metric and measure,} or an \index{mmspace@$mm$-space}
{\em $mm$-space}, is a triple,
$(X,d,\mu)$, consisting of a set $X$, a metric $d$ on $X$,
and a probability Borel measure $\mu$ on the metric space $(X,d)$.
\end{definition}

For a subset $A$ of a metric space $X$ and an $\e>0$, denote by $A_\e$ the
$\e$-neighbourhood of $A$ in $X$.

\begin{definition}[{\it ibid.}]
\label{def:levy}
A family ${\mathcal X}=(X_n,d_n,\mu_n)_{n\in\N}$ 
of $mm$-spaces is called a {\em L\'evy family} if,
whenever Borel subsets $A_n\subseteq X_n$ satisfy 
\[\liminf_{n\to\infty}\mu_n(A_n)>0,\]
one has for every $\e>0$
\[
\lim_{n\to\infty}\mu_n((A_n)_\e) = 1.
\label{condi}
\]
\end{definition}

The concept of a L\'evy family can be reformulated in many
equivalent ways. For example, it is not difficult to see
that a family ${\mathcal X}$ as above is L\'evy if and only if
for every $\e>0$, whenever $A_n,B_n$ are Borel subsets of $X_n$ satisfying
\[\mu_n(A_n)\geq\e,~~\mu_n(B_n)\geq \e,\]
one has $d(A_n,B_n)\to 0$ as $n\to\infty$.

This is formalized using the notion of \textit{separation distance,} proposed
by Gromov (\cite{Gr}, Section 3$\frac 12$.30). Given numbers 
$\kappa_0,\kappa_1,\ldots,\kappa_N>0$, one defines the invariant
\[{\mathrm{Sep}}\,(X;\kappa_0,\kappa_1,\ldots,\kappa_N)\]
as the supremum of all $\delta$ such that $X$ contains Borel subsets
$X_i$, $i=0,1,\ldots,N$ with $\mu(X_i)\geq\kappa_i$, every two of which are
at a distance $\geq\delta$ from each other.
Now a family ${\mathcal X}=(X_n,d_n,\mu_n)_{n\in\N}$ is a L\'evy family if
and only if for every $0<\e<\frac 12$, one has
\[{\mathrm{Sep}}\,(X_n;\e,\e)\to 0\mbox{ as }n\to\infty.\]
The reader should consult Ch. 3$\frac 12$ in \cite{Gr} for numerous other
characterisations of L\'evy familes of $mm$-spaces.

We will state just one more such reformulation. 
It is an easy exercise to show that 
in the Definition \ref{def:levy} of a L\'evy family
it is enough to
assume that the values $\mu_n(A_n)$ are bounded away from zero by
$1/2$ (or by any other fixed constant strictly between zero and one). 
In other words, 
a family ${\mathcal X}$ is a L\'evy family if
and only if, whenever Borel subsets $A_n\subseteq X_n$ 
satisfy $\mu_n(A_n)\geq 1/2$, one has for every $\e>0$
\[
\lim_{n\to\infty}\mu_n(A_n)_\e = 1.
\label{condi2}
\]

This leads to the following concept \cite{M86,MS}, 
providing convenient quantitative bounds on the rate of convergence of 
$\mu_n(A_n)_\e$ to one. 

\begin{definition}
Let $(X,d,\mu)$ be a space with metric and measure.
The {\em concentration function} of $X$,
denoted by $\alpha_X(\e)$, is a
real-valued function on the positive axis $\R_+=[0,\infty)$,
defined by letting $\alpha(0)=1/2$ and for all $\e>0$
\[\alpha_X(\e)=
1-\inf\left\{\mu\left(B_\e\right) \colon
B\subseteq X, ~~ \mu(B)\geq\frac{1}{2}\right\}.
\]
\label{defconcfn}
\end{definition}

Thus, a family ${\mathcal X}=(X_n,d_n,\mu_n)_{n\in\N}$ 
of $mm$-spaces is a L\'evy family if and only if 
\[\alpha_{X_n}\to 0\mbox{ pointwise on $(0,+\infty)$ as }n\to\infty.\]
A L\'evy family is called {\em normal} 
if for suitable constants $C_1,C_2>0$,
\[\alpha_{X_n}(\e)\leq C_1e^{-C_2\e^2n}.\]

\begin{example}
The Euclidean spheres $\s^n$, $n\in\N_+$ of unit
radius, equipped with the Haar measure (translation-invariant probability
measure) and Euclidean (or geodesic) distance, form a normal L\'evy family. 
\end{example}

\begin{definition}[Gromov and Milman \cite{GrM}]
A metrizable 
topological group $G$ is called a \emph{L\'evy group} if 
it contains an increasing chain of compact subgroups
\[G_1<G_2<\ldots<G_n<\ldots,\]
having an everywhere dense union in $G$ and such that for some right-invariant
compatible metric $d$ on $G$ the groups $G_n$, equipped with the normalized
Haar measures and the restrictions of the metric $d$, form a L\'evy
family. 
\end{definition}

The above concept admits a number of generalizations, in particular it makes 
perfect sense for non-metrizable, non-separable topological groups as well. 
In fact, in the definintion of a L\'evy group it
is the uniform structure on $G$ that matters rather than a metric. Namely,
one can easily prove the following.

\begin{proposition}
\label{p:va}
Let $G$ be a metrizable topological group containing an increasing chain of
compact subgroups $(G_n)$ with everywhere dense union. The subgroups
$(G_n)$ form a L\'evy family with regard to the normalized
Haar measures and the restrictions of some right-invariant metric $d$ on $G$
if and only if for every neighbourhood of identity, $V$, in $G$ and every
collection of Borel subsets $A_n\subseteq G_n$ with the property
$\mu_n(A_n)\geq 1/2$ one has
\[\lim_{n\to\infty}\mu_n(VA_n)=1.\]
\qed
\end{proposition}

Examples of presently known L\'evy groups can be found in 
\cite{GrM,MS,Gl1,P99,L,GP2,GTW}.

The following result had been also established in \cite{GrM}, and one can give
numerous alternative proofs to it, cf. e.g. \cite{Gl1,P99,GP2,P05b}.

\begin{theorem} Every L\'evy group is extremely amenable. \qed
\end{theorem}

The concept of a L\'evy group is stronger than that of an
extremely amenable group. Typically, examples of extremely amenable groups
coming from combinatorics as groups of automorphisms of infinite Fra\"\i ss\'e
order structures \cite{KPT} are not L\'evy groups, because they contain no
compact subgroups whatsoever. Even the dynamical behaviour of L\'evy groups
has been shown by Glasner, Tsirelson and Weiss \cite{GTW}
to differ considerably from
that of the rest of extremely amenable groups. 

The main theorem of this
article (Th. \ref{th:isoulevy}) states
that the group $\Iso(\Ur)$ is a L\'evy group rather than merely an extremely
amenable one.

The monograph \cite{P05} by the present author
provides an introduction to the theory of
extremely amenable groups and its links with geometric functional analysis
and combinatorics. The Urysohn metric space can also be found there. In fact,
the book also contains Theorem \ref{th:isoulevy} (cf. Section 3.4.3). 
The reason to
have this theorem published in the present Proceedings is twofold: firstly, the book \cite{P05}, published in Brazil, is not readily available, and secondly, the proof of Theorem \ref{th:isoulevy} as presented there is not sufficiently accurate.

\section{Approximating $\Iso(\Ur)$ with finite subgroups\label{s:lfs}}

Let $\Gamma=(V,E)$ be an (undirected, simple) graph, where $V$ is the set of
vertices and $E$ is the set of edges. 
A {\em weight} on $\Gamma$ is an assignment of a non-negative real number to
every edge, that is, a function $w\colon E\to\R_+$. The pair $(\Gamma,w)$ forms
a \index{graph!weighted}
{\em weighted graph}. The \index{path pseudometric}
{\em path pseudometric} on a connected weighted 
graph $(\Gamma,w)$ is the maximal pseudometric on 
$\Gamma$ with the property $d(x,y)=w(x,y)$ for any pair of adjacent vertices $x,y$.
Equivalently, the value of $\rho(x,y)$ is
given for each $x,y\in V$ by
\begin{equation}
\rho(x,y) = \inf\sum_{i=0}^{N-1} d(a_i,a_{i+1}),
\label{eq:minimum}
\end{equation}
where the infimum is taken over all positive natural $N$ and all finite sequences
of vertices
$x=a_0, a_1,\ldots,a_{N-1},a_N=b$, with the property that $a_i$ and $a_{i+1}$
are adjacent for all $i$. Notice that here we allow for sequences of length
one, in which case the sum above is empty and returns value zero, the distance
from a vertex to itself.


In particular, if every edge is assigned the weight one, the corresponding
path pseudometric is a metric, called the {\em path metric} on $\Gamma$.

Let $G$ be a group, and $V$ a generating subset of
$G$. Assume that $V$ is symmetric ($V=V^{-1}$) and contains the identity. The 
{\em Cayley graph} associated to the pair $(G,V)$ has all elements of the group
$G$ as vertices, and two of them, $x,y\in G$, $x\neq y$, are adjacent
if and only if $x^{-1}y\in V$. The Cayley graph is connected.
The corresponding path metric on $G$ is called the \index{word distance}
{\em word distance} with regard to the generating set $V$. 

If $V$ is an arbitrary generating subset of $G$, then the word distance with
regard to $V$ is defined as that with regard to $V\cup V^{-1}\cup\{e\}$.
The value of the word distance between $e$ and an element $x$ is
called the \index{reduced length}
{\em reduced length} of $x$ with regard to the generating set $V$,
and denoted $\ell_V(x)$. It is simply the smallest integer $n$ such that
$x$ can be written as a product of $\leq n$ elements of $V$ and their
inverses. Since the identity $e$ of the group $G$ is represented,
as usual, by an empty word, one has $V^0=\{e\}$ and $\ell_V(e)=0$.

\begin{lemma} 
\label{l:maxim}
Let $G$ be a group equipped with a left-invariant pseudometric,
$d$. Let $V$ be a finite generating 
subset of $G$ containing the identity.
Then there is the maximal pseudometric, $\rho$,
among all left-invariant pseudometrics on $G$, whose restriction to $V$
is majorized by $d$. The restrictions of $\rho$ and $d$ to $V$ coincide.
If $d\vert_V$ is a metric on $V$, then $\rho$ is a metric as well, and
for every $\e>0$ there is an $N\in\N$ such that  
$\ell_V(x)\geq N$ implies $\rho(e,x)\geq\e$.
\end{lemma}

\begin{pf} 
Make the Cayley graph $\Gamma$ associated to the pair $(G,V^{-1}V)$ into a
weighted graph, by assigning to every edge $(x,y)$, $x^{-1}y\in V^{-1}V$,
the value $d(x,y)\equiv d(x^{-1}y,e)$.
Denote by $\rho$ the corresponding path pseudometric on the weighted
graph $\Gamma$. To prove the left-invariance of $\rho$, let $x,y,z\in G$.
Consider any sequence of elements of $G$,
\begin{equation}
x=a_0, a_1,\ldots,a_{N-1},a_N=y,
\label{eq:seq}
\end{equation}
where $N\in\N$ and $a_i^{-1}a_{i+1}\in V^{-1}V$, $i=0,1,\ldots,n-1$.
Since for all $i$ the elements $za_i,za_{i+1}$ are adjacent in the Cayley
graph ($(za_i)^{-1}za_{i+1}=a_i^{-1}a_{i+1}\in V^{-1}V$), one has
\begin{eqnarray*}
d(zx,zy) &\leq & 
\sum_{i=0}^{n-1} d(za_i,za_{i+1}) \\
&=& \sum_{i=0}^{n-1} d(a_i,a_{i+1}),
\end{eqnarray*}
and taking the infimum over all sequences as in Eq. (\ref{eq:seq}) on both 
sides, one concludes $d(zx,zy)\leq d(x,y)$, which of course implies the
equality.

For every $x,y\in V$ one has $x^{-1}y\in V^{-1}V$ and consequently
$\rho(x,y)=\rho(x^{-1}y,e)\leq d(x^{-1}y,e)=d(x,y)$.
Now let $\varsigma$ be any left-invariant pseudometric on $G$
whose restriction to $V$ is majorized by $d$. 
If $a,b\in G$ are such that $a^{-1}b\in V^{-1}V$, then for some $c,d\in V$
one has $a^{-1}b=c^{-1}d$, and
\begin{eqnarray*}
\varsigma(a,b) = \varsigma(a^{-1}b,e)
= \varsigma(c^{-1}d,e) 
=
\varsigma(c,d) 
&\leq& d(c,d)  = d(c^{-1}d,e)  =d(a^{-1}b,e) = d(a,b).
\end{eqnarray*}
For every sequence as in 
Eq. (\ref{eq:seq}), one now has
\begin{eqnarray*}
\varsigma(x,y) \leq  \sum_{i=0}^{n-1} \varsigma(a_i,a_{i+1}) 
\leq  \sum_{i=0}^{n-1} d(a_i,a_{i+1}),
\end{eqnarray*}
and by taking the infimum over all such finite sequences on both sides, one
concludes  
\[\varsigma(x,y)\leq \rho(x,y),\]
that is, $\rho$ is maximal among all left-invariant pseudometrics whose
restriction to $V$ is majorized by $d$. In particular, $\rho\geq d$, which
implies $\rho\vert_V=d\vert_V$. 

Assuming that $d\vert_V$ is a metric, all the weights on the Cayley graph
$\Gamma$ as above
assume strictly positive values, and consequently $\rho$ is a metric.
As we have already noticed,
for every $x,y\in G$ with the property $x^{-1}y\in V^{-1}V$, the 
value $d(x,y)$ is of the form $d(a,b)$ for suitable $a,b\in V$. Consequently,
there exists the smallest value taken by $d$ between
pairs of distinct elements $x,y\in G$ with the property $x^{-1}y\in V^{-1}V$,
and it is strictly positive. Denote this value by $\delta$. 
Clearly, for every $x\in G$ one has
$\rho(e,x)\geq\delta \ell_{V^{-1}V}(x)\geq (\delta/2)\ell_V(x)$, and the proof is finished.
\qed\end{pf}

Next we are going to get rid of the restrictions on $V$. The price to pay is
to agree that all pseudometrics will be bounded by $1$.
In the following lemma, $\ell_V(x)$ will denote the word length of $x$
with regard to $V$ if $x$ is contained in the subgroup generated by $V$,
and $\infty$ otherwise.

\begin{lemma} 
\label{l:maxim2}
Let $G$ be a group equipped with a left-invariant pseudometric,
$d$, whose values are bounded by $1$. Let $V$ be a finite  
subset of $G$.
Then there is the maximal pseudometric, $\rho$,
among all left-invariant pseudometrics on $G$, bounded by one and
whose restriction to $V$
is majorized by $d$. The restrictions of $\rho$ and $d$ to $V$ coincide.
If $d\vert_V$ is a metric on $V$, then $\rho$ is a metric on $G$.
\end{lemma}

\begin{pf} 
The set $\Psi$ of all left-invariant pseudometrics on $G$ bounded by one
and whose restrictions to $V$ are majorized by $d$ is non-empty ($d\in\Psi$),
and contains the maximal element, $\rho$, given by 
$\rho(x,y)=\sup_{\varsigma\in\Psi}\varsigma(x,y)$. Obviously,
$\rho\vert V =d\vert V$.
To verify the last assertion, let $\delta$ be the smallest strictly positive
value of the form $d(x,y)$, $x,y\in V$, $x\neq y$. Let $\varsigma$ now denote
the metric on $G$ taking values $0$ and $\delta$. Denote by 
$\langle V\rangle$ the subgroup of $G$ generated by $V$. 
According to Lemma \ref{l:maxim},
there exists the maximal metric $\varsigma_1$ on the subgroup $\langle V\rangle$ 
of $G$ generated by $V$ whose restriction to 
$V\cup V^{-1}\cup\{e\}$
only takes the values $0$ or $\delta$. Define a pseudometric $\varsigma_2$
on all of $G$ by the rule
\[\varsigma_2(x,y) = \left\{\begin{array}{ll} \min\{1, \varsigma_1(x,y)\},& 
\mbox{ if }x^{-1}y\in \langle V\rangle, \\
1, &\mbox{otherwise.}
\end{array}\right.
\]
Since $\varsigma_1\vert V \leq d\vert V$,
it follows that $\varsigma_2\vert V \leq\rho$, and thus $\rho$ is a metric.
\qed\end{pf}

\begin{lemma} 
\label{l:factorm}
Let $\rho$ be a left-invariant pseudometric on a group $G$, and let
$H\triangleleft G$ be a normal subgroup. The formula
\begin{eqnarray}
\label{infimum}
\bar \rho(xH,yH) &:=& \inf_{h_1,h_2\in H} \rho(xh_1,yh_2) \\
&\equiv& \inf_{h_1,h_2\in H} \rho(h_1x,h_2y) \nonumber\\
&\equiv& \inf_{h\in H} \rho(hx,y)  \nonumber
\end{eqnarray}
defines a left-invariant pseudometric on the factor-group
$G/H$. This is the largest pseudometric on $G/H$ with respect to which the
quotient homomorphism $G\to G/H$ is $1$-Lipschitz.
\end{lemma}

\begin{pf}
The triangle inequality follows from the fact that,
for all $h'\in H$,
\begin{eqnarray*}
\bar \rho(xH,yH) &=&
\inf_{h\in H} \rho(hx,y) \\
&\leq &
\inf_{h\in H} [ \rho(hx,h'z) +  \rho(h'z,y) ]\\
&=& \inf_{h\in H} \rho(hx,h'z) +  \rho(h'z,y) \\
&=& \inf_{h\in H} \rho(h'^{-1}hx,z) +  \rho(h'z,y) \\
&=& \bar\rho(xH,zH) + \rho(h'z,y),
\end{eqnarray*}
and the infimum of the r.h.s. taken over all  $h'\in H$ equals
$\bar\rho(xH,zH) + \bar\rho(zH,yH)$. Left-invariance of $\bar\rho$ is obvious.
If $d$ is a pseudometric on $G/H$ making the quotient homomorphism into a
$1$-Lipschitz map, then $d(xH,yH)\leq \rho(xh_1,yh_2)$ for all $x,y\in G$,
$h_1,h_2\in H$, and therefore $d(xH,yH)\leq\bar\rho(xH,yH)$.
\qed\end{pf}

We will make a distinction between the notion of a
{\em distance-preserving map} $f\colon X\to Y$ 
between two pseudometric spaces, which has the property $d_Y(fx,fy)=d_X(x,y)$
for all $x,y\in X$, and an {\em isometry}, that is, a distance-preserving
bijection. 

Let $G$ be a group.
For every left-invariant bounded pseudometric $d$ on $G$, denote
$H_d=\{x\in G\colon d(x,e)=0\}$, and let $\hat d$ be the metric on
the left coset space $G/H_d$ given by $\hat d(xH_d,yH_d)=d(x,y)$.
The metric $\hat d$ is invariant under left translations by elements of $G$.
We will denote the metric space $(G/H_d,\hat d)$, equipped with the left
action of $G$ by isometries, simply by $G/d$.

A distance-preserving map need not be an isometry. For instance,
if $d$ is a left-invariant pseudometric on a group $G$, then the natural
map $G\to G/d$ is distance-preserving, onto, but not necessarily an injection.

A group $G$ is {\em residually finite} if it admits a separating
family of homomorphisms into finite groups, or, equivalently, if for every
$x\in G$, $x\neq e$, there exists a normal subgroup $H\triangleleft G$ of
finite index such that $x\notin H$. Every free group is residually
finite (cf. e.g. \cite{MKS}), and the free product of two residually finite groups is residually finite \cite{Gru}.

\begin{lemma} 
\label{l:rfg}
Let $G$ be a residually finite group equipped with a
left-invariant pseudometric $d\leq 1$, and let $V\subseteq G$ be a finite 
subset.
Suppose the restriction $d\vert_V$ is a metric, and
let $\rho$ be the maximal left-invariant metric on $G$
bounded by one with $\rho\vert_V=d\vert_V$.
Then there exists a normal subgroup $H\triangleleft G$ of finite index
with the property that the restriction of the quotient homomorphism 
$G\to G/H$ to $V$ is an isometry with regard to $\rho$ and the 
quotient pseudometric $\bar\rho$ (which is in fact a metric).
\end{lemma}

\begin{pf} 
Let $\delta>0$ be the smallest distance between any pair of distinct elements
of $V$.
Let $N\in\N_+$ be even and such that $(N-2)\delta\geq 1$. 
The subset formed by all words of length $\leq N$ in $V$ is
finite, and, since the intersection of finitely many 
subgroups of finite index has finite
index (Poincar\'e's theorem), one can choose a normal subgroup 
$H\triangleleft G$ of finite index containing no words of $V$-length $\leq N$
other than $e$. 
Let $x,y\in V$ and $h\in H$, $h\neq e$. If 
$y^{-1}hx\notin\langle V\rangle$, then 
$\rho(hx,y)=1$. If $y^{-1}hx\in\langle V\rangle$, then 
the reduced $V$-length of the word $y^{-1}hx$ is
$\geq N-2$, and consequently
\[\rho(hx,y)=\rho(y^{-1}hx,e)\geq \min\{1,(N-2)\delta\}\geq 1.\]
(Otherwise there would exist a representation 
\[y^{-1}hx = v_1v_2\ldots v_{k}\]
with $v_i\in V$ and 
$d(e,v_1)+\sum_{i=1}^{k-1} d(v_i,v_{i+1})< (N-2)\delta$, that is,
$k< N-2$.)

In either case, the distance $\bar\rho(xH,yH)$ between cosets is realized on
the representatives $x,y$: 
\[\bar\rho(xH,yH)=\rho(x,y).\]
The factor-pseudometric $\bar\rho$ on $G/H$ is, according to 
Lemma \ref{l:factorm}, the largest pseudometric making the factor-map $\pi$ 
1-Lipschitz. We claim that $\bar\rho$ is the largest left-invariant 
pseudometric on $G/H$, bounded by one, whose restriction to $V$ coincides with
the metric on $V$. Indeed, denoting such a pseudometric by $\varsigma$, one sees
that $\varsigma\circ\pi$ is a left-invariant pseudometric on $G$, bounded by
one, and whose restriction to $V$ equals $d_\xi\vert V$. It follows that
$\varsigma\circ\pi\leq\rho$, thence $\varsigma\leq\bar\rho$ and the two
coincide. Now Lemma \ref{l:maxim2} tells us that $\bar\rho$ is a metric.
\qed\end{pf}

The following concept, along with the two subsequent results, forms a
powerful tool in the theory of the Urysohn space.

\begin{definition}
[Uspenskij \cite{Usp98}]
One says that a metric subspace $Y$ is {\rm $g$-embedded} into a metric space 
$X$ if there exists an embedding of topological groups
$e\colon\Iso(Y)\hookrightarrow\Iso(X)$ with the property that
for every $h\in\Iso(Y)$ the isometry $e(h)\colon X\to X$ is
an extension of $h$:
\[e(h)\vert_X = h.\]
\end{definition}

\begin{proposition}
[Uspenskij \cite{Usp90,Usp98}]
Each separable metric space $X$ admits a $g$-embedding into 
the complete separable Urysohn metric space $\Ur$.
\label{p:sg-emb}
\qed
\end{proposition}

Every isometry between two compact metric subspaces of the Urysohn space
$\Ur$ extends to a global self-isometry of $\Ur$ (it was first established
in \cite{hui}). 
Together with Proposition \ref{p:sg-emb}, this fact
immediately leads to the following result.

\begin{proposition}
\label{c:conven}
Each isometric embedding of a compact metric space into $\Ur$ is a
$g$-embedding. \qed
\end{proposition}

Recall that an action of a group $G$ on a set $X$ is \textit{free} if
for all $g\in G$, $g\neq e$ and all $x\in X$, one has $g\cdot x\neq x$.
Here comes the main technical result of this paper.

\begin{lemma}
Let $X$ be a finite subset of the Urysohn space $\Ur$, and let a finite
group $G$ act on $X$ freely by isometries. 
Let $f$ be an isometry of $\Ur$,
and let $\e>0$. There exist a finite group $\tilde G$ containing $G$ as
a subgroup, an element $\tilde f\in \tilde G$, and a finite 
metric space $Y$, $X\subseteq Y\subset\Ur$, upon which 
$\tilde G$ acts freely by isometries, extending the 
original action of $G$ on $X$ and so that for all $x\in X$ one has
$d(\tilde f x,fx)<\e$.
\label{l:approx}
\end{lemma}

\begin{pf} 
Without loss in generality, one can assume that the image $f(X)$ does not
meet $X$, by replacing $f$, if necessary,
with an isometry $f^\prime$ such that the image $f^\prime(X)$ does not
intersect $X$, and yet for every $x\in X$ one has $d_{\Ur}(f(x),f^\prime(x))<\e$.
By renormalizing the distance if necessary, we will further assume that
the diameter of the set $X\cup f(X)$ does not exceed $1$. 

Since every compact subset of $\Ur$ 
such as $X$ is $g$-embedded into the Urysohn space
(Proposition \ref{c:conven}),
one can extend the action of $G$ by isometries from $X$ to all of $\Ur$.

Choose any element $\xi\in \Ur$ at a distance $1$ from every element of
$X\cup f(X)$. Let $\Theta=X/G$ denote the set of $G$-orbits of $X$. For 
each $\theta\in\Theta$, choose an element $x_\theta\in\theta$ and an
isometry $f_\theta$ of $\Ur$ in such a way that $f_\theta(\xi)=x_\theta$. 
Let $n=\abs\Theta$, and let $F_n$ be the free group on $n$ generators
which we will denote likewise $f_\theta$, $\theta\in\Theta$. 

Finally, denote by $f$ a generator of the group $\Z$, and let
$F=G\ast F_n\ast\Z$ be the free product of three groups.

There is a unique homomorhism $F\to \Iso(\Ur)$, which sends all elements
of $G\cup\{f_\theta\colon\theta\in\Theta\}\cup\{f\}$ to the corresponding 
self-isometries of $\Ur$. In this way, $F$ acts on $\Ur$ by isometries.
Denote 
\[V=\{g\circ f_\theta\colon g\in G,~\theta\in\Theta\}
\cup \{f\circ g\circ f_\theta\colon g\in G,~\theta\in\Theta\}.\]
The formula
\[d_\xi(g,h):= \min\{1,d_\Ur(g(\xi),h(\xi))\},~~g,h\in F,\]
defines a left-invariant pseudometric $d_\xi$ on the group $F$, bounded by $1$.

Denote by $\ev\colon F\to\Ur$ the evaluation map $\phi\mapsto \phi(\xi)$.
The restriction $\ev\vert V$ is
an isometry between $V$, equipped with the restriction of the 
pseudometric $d_\xi$,
and $X\cup f(X)$. Also notice that the restriction 
$\ev\vert \{g\circ f_\theta\colon g\in G,~\theta\in\Theta\}$ establishes an
isomorphism of $G$-spaces between the latter set (upon which
$G$ acts by left multipication in the group $F$) and $X$. Both properties
take into account the freeness of the action of $G$ on $X$.

The restriction of the pseudometric $d_\xi$ to $V$ is a metric.
Let $\rho$ be the maximal left-invariant metric on $F$ bounded
by $1$ such that $\rho\vert_V=d_{\xi}\vert_V$. (Lemma \ref{l:maxim2}.)

The group $F$, being the free product of three residually finite groups,
is residually finite, and so we are under the assumptions of Lemma \ref{l:rfg}.
Choose a normal subgroup $H\triangleleft F$ of finite index in such a way
that if the finite factor-group $F/H$ is equipped with the
factor-pseudometric $\bar\rho$, then the restriction of the factor-homomorphism
$\pi\colon F\to F/H$ to $V$ is an isometry. This $\bar\rho$ is then a metric.
In addition, by replacing $H$ with a smaller normal subgroup of finite index
if necessary, one can clearly choose
$H$ so that $H\cap G=\{e\}$, and thus $\pi\vert G$ is a monomorphism.

The finite group $\tilde G=F/H$ acts on itself by left translations, 
and this action is a free action by isometries on the finite metric space
$Y = (F/H,\bar\rho)$. The metric space $X\cup f(X)$ embeds into $Y$ as a metric
subspace through the isometry $\pi\circ \ev$, and $\tilde f\vert X = f\vert X$.
Finally, $G$ is a subgroup of $\tilde G$, and $X$ is contained inside $Y$
as a $G$-space. 

Finally, the embedding of $X\cup f(X)$ (considered as a subspace of $Y$) can
be extended over $Y$, so we can view $Y$ as a metric subspace of $\Ur$, and
it is a $g$-embedding by Proposition \ref{c:conven}.
\qed\end{pf}

Now we are ready to give an alternative proof of the following result of
Vershik. Recall that a group $G$ is \textit{locally finite} if every
finitely generated subgroup of $G$ is finite. A countable group is
locally finite if and only if it is the union of an increasing chain of
finite subgroups.

\begin{theorem}[Vershik \cite{Ver05b}]
\label{th:lf}
The isometry group $\Iso(\Ur)$ of the Urysohn space,
equipped with the standard Polish topology, contains an everywhere
dense locally finite countable subgroup. 
\end{theorem}

\begin{pf} Choose an everywhere dense subset $F=\{f_i\colon i\in\N_+\}$
of $\Iso(\Ur)$ and a point $x_1\in\Ur$. 

Let $G_1=\{e\}$ be a trivial group, trivially acting on $\Ur$ by isometries.
Clearly, the restriction of this action on the $G_1$-orbit of $\{x_1\}$
is free.

Assume that for an $n\in\N$ one has chosen recursively
a finite group $G_n$, an
action $\sigma_n$ by isometries on $\Ur$, and a collection of points
$\{x_1,\ldots,x_{2^n}\}$
in such a way that the restriction of the action $\sigma_n$ to the $G_n$-orbit of
$\{x_1,x_2,\ldots,x_{2^n})$ is free.

Using Lemma \ref{l:approx}, choose a finite group 
$G_{n+1}$ containing (an isomorphic copy) of $G_n$, 
an element $\tilde f_n\in G_{n+1}$ 
and an action $\sigma_{n+1}$ of $G_{n+1}$ on $\Ur$ by isometries such that 
for every $j=1,2,\ldots,2^n$ and each $g\in G_n$ one has 
\[\sigma_n(g)x_j=\sigma_{n+1}(g)x_j,\]
the elements $x_{2^n+j}=\tilde f_n(x_j)$, $j=1,2,\ldots,2^n$ are all
distinct from any of $x_i$, $i\leq 2^n$,
the restriction of the action of $G_{n+1}$ on the $G_{n+1}$-orbit of
$\{x_0,x_1,\ldots,x_{2^{n+1}}\}$ is free, and
\[d_{\Ur}(f_n(x_j),\tilde f_n(x_j))<2^{-n},~~j\leq 2^n.\]

The subset $X=\{x_i\colon i\in\N_+\}$ is everywhere dense in
$\Ur$. Indeed, for each $n\in\N$
the subset $\{f_i(x_n)\colon i\geq n\}$ is everywhere dense in $\Ur$,
and since it is contained in the $2^{-n}$-neighbourhood of
$\{\tilde f_i(x_n)\colon i\geq n\}\subset X$, the statement follows.

The group $G=\cup_{i=1}^\infty G_n$ is locally finite. Now let $g\in G$. For
every $i\in\N_+$, the value $g\cdot x_i$ is well-defined as the limit
of an eventually constant sequence, and determines an
isometry from an everywhere dense subset $X\subset\Ur$ into $\Ur$.
Consequently, it extends uniquely to an isometry from $\Ur$ into itself.
If $g,h\in G$, then the isometry determined by $gh$ is
the composition of isometries determined by $g$ and $h$: every $x\in X$
has the property $(gh)(x) = g(h(x))$, once $x=x_i$, $i\leq N$, and
$g,h\in G_N$, and this property extends over all of $\Ur$. Thus,
$G$ acts on $\Ur$ by isometries (which are therefore onto). 

Finally, notice that $G$ is everywhere dense in $\Iso(\Ur)$. It is enough
to consider the basic open sets of the form
\[\{f\in\Iso(\Ur)\colon d(f(x_i),g(x_i))<\e,~~i=1,2,\ldots,n\},\]
where $g\in\Iso(\Ur)$, $n\in\N$, and $\e>0$.
Since $F$ is everywhere dense in $\Iso(\Ur)$, there is an $m\in\N$
with $n\leq 2^{m-1}$, $2^{-m}<\e/2$, and
$d(f_m(x_i),g(x_i))<\e/2$ for all $i=1,2,\ldots,n$. One concludes:
$d(\tilde f_m(x_i),g(x_i))<\e$ for $i=1,2,\ldots,n$, and 
$\tilde f_m\in G_m\subset G$, which settles the claim.
\qed\end{pf}

A further refinement of our argument leads to another approximation
theorem \ref{th:isoulevy}, which states that $\Iso(\Ur)$ is a L\'evy group and
forms the central result of the present paper. The proof will 
interlace the recursion steps in the proof of
Theorem \ref{th:lf} with an adaptation of an idea used in the proof of the
following result to obtain, 
historically, the second ever example of a L\'evy
group, after $U(\ell^2)$.

\begin{theorem}[Glasner \cite{Gl1}; Furstenberg and Weiss (unpublished)]
Let $G$ be a compact metric group, and let $d$ be an invariant metric on $G$. 
The group $L^1([0,1];G)$ of all equivalence classes of
Borel maps from the unit interval $[0,1]$ to $G$, equipped with the
metric $d_1(f,g) =\int_0^1 d(f(x),g(x))dx$, is a L\'evy group.
\qed
\end{theorem}

The following well-known and important result is being established using the
probabilistic techniques (martingales). 
(Cf. the more general Theorem 7.8 in \cite{MS} or Theorem 4.2 in \cite{L}.)

\begin{theorem}
Let $(X_i,d_i,\mu_i)$, $i=1,2,\ldots,n$ be metric spaces with
measure, each having finite diameter $a_i$.
Equip the product $X(i)=\prod_{i=1}^n X_i$ with the product measure
$\otimes_{i=1}^n\mu_i$ and the $\ell_1$-type (Hamming) metric
\[d(x,y) = \sum_{i=1}^n d_i(x_i,y_i).\]
Then the concentration function of $X$ satisfies
\[
\alpha_X(\e)\leq  2 e^{-\e^2/8\sum_{i=1}^n a_i^2}.
\]
\label{th:product}
\qed
\end{theorem}

Let us consider the following particular case.
Let $(X,d)$ be a finite metric space, and let $Z$ be a finite set
equipped with the normalized counting measure $\mu_{\sharp}$, that is,
$\mu_{\sharp}(A) = \abs A/\abs Z$. We will equip the collection
$X^Z$ of all maps from $Z$ to $X$ with the $L_1(\mu_{\sharp})$-metric:
\[d_1(f,g) = \int_Z d(f(z),g(z))~d\mu_{\sharp}(z).\]
This is just the $\ell_1$-metric normalized:
\[d_1(f,g) = \frac{1}{\abs Z}\sum_{z\in Z} d(f(z),g(z)).\]
It is also known as the (generalized) {\em normalized Hamming distance.}
In particular, if $a=\diam(Z)$ is the diameter of $Z$, then
the diameter of every ``factor'' of the form
$\{z\}\times Z$ is $a/n$, and Theorem \ref{th:product} gives the following.

\begin{corollary}
Let $(X,d)$ is a finite metric space of diameter $a$ and let $n\in\N$.
Let the metric space 
$X^n$ be equipped with the normalized counting measure and the
normalized Hamming distance.
Then the concentration function of $mm$-space $X^n$ satisfies
\[
\alpha_{X^n}(\e)\leq  2 e^{-n\e^2/8a^2}.
\]
\label{c:xtoz}
\qed
\end{corollary}

Notice that $X^n$ with the above metric contains an isometric copy of $X$,
consisting of all constant functions. 

If a finite group $G$ acts on a finite
metric space $X$ by isometries, then this action naturally extends to an
action of $G^n$ on $X^n$ by isometries, where the latter set is equipped with
the normalized Hamming, or $L_1(\mu_{\sharp})$, metric. If the action of $G$
on $X$ is free, then so is the action of $G^n$ on $X^n$.

\begin{theorem}
\label{th:isoulevy}
The isometry group $\Iso(\Ur)$ of the Urysohn space,
equipped with the standard Polish topology, is a L\'evy group.
Moreover, the groups in the approximating L\'evy family can be chosen finite.
\end{theorem} 

\begin{pf}
As in the proof of Theorem \ref{th:lf},
choose an everywhere dense subset $F=\{f_i\colon i\in\N_+\}$
of $\Iso(\Ur)$ and a point $x_1\in\Ur$. Set $G_1=\{e\}$ and $X_1=\{x_1\}$.
Assume that for an $n\in\N_+$
a finite group $G_n$, an action $\sigma_n$ by isometries on $\Ur$, and
a finite $G_n$-invariant subset $X_n\subset\Ur$ have been chosen. Also assume
that $G_n$ acts on $X_n$ freely. 
Let $a_n$ be the diameter of $X_n$. Choose $m_n\in\N$ so that
\begin{equation}
m_n\geq 8 a_n^2n.
\label{eq:growth}
\end{equation}
The finite metric space $\tilde X_n =X^{m_n}_n$ (with the $L_1(\mu_{\sharp})$-metric)
contains $X_n$ as a subspace of constant functions,
therefore one can embed $\tilde X_n$
into $\Ur$ so as to extend the embedding $X_n\hookrightarrow\Ur$ (the finite
injectivity of $\Ur$).

The group $\tilde G_n = G^{m_n}_n$ acts on the metric space $\tilde X_n$ freely
by isometries.
Since every embedding of a
compact subspace into $\Ur$ is a $g$-embedding, one can simultaneously
extend the action of $\tilde G_n$ to a global action, $\tilde\sigma_n$,
on $\Ur$ by isometries.
Now construct the group $G_{n+1}$ and its action $\sigma_{n+1}$
by isometries exactly as in the proof of Theorem \ref{th:lf}, but beginning
with $\tilde G_n$ instead of $G_n$ and $\tilde X_n$ instead of 
$\{x_1,\ldots,x_{2^n}\}$. Namely, using 
Lemma \ref{l:approx}, choose a finite group 
$G_{n+1}$ containing (an isomorphic copy) of $\tilde G_n$, 
an element $\tilde f_n\in G_{n+1}$ 
and an action $\sigma_{n+1}$ of $G_{n+1}$ on $\Ur$ by isometries such that 
for every $x\in \tilde X_n$ and each $g\in G_n$ one has 
\[\sigma_n(g)x=\sigma_{n+1}(g)x,\]
the sets $\tilde f_n(\tilde X_n)$ and $\tilde X_n$ are disjoint,
the restriction of the action of $G_{n+1}$ on the $G_{n+1}$-orbit of
$\tilde X_n$ is free, and
\[d_{\Ur}(f_n(x),\tilde f_n(x))<2^{-n}~~\mbox{for all $x\in\tilde X_n$}.\]
Denote $X_{n+1}= G_{n+1}\cdot \tilde X_n$. The step of recursion is accomplished.

The union 
$G=\cup_{i=1}^\infty G_n=\cup_{i=1}^\infty\tilde G_n$ is, like
in the proof of Theorem \ref{th:lf},
an everywhere dense locally finite subgroup of $\Iso(\Ur)$, and it only 
remains to show that the groups $\tilde G_n$, $n\in\N_+$, form a L\'evy
family
with regard to the uniform structure inherited from $\Iso(\Ur)$.

First, consider the groups 
$\tilde G_n=G^{m_n}_n$ equipped with the $L(\mu_{\sharp})$-metric formed
with regard to the \textit{discrete} (that is, $\{0,1\}$-valued)
metric on $G_n$. 
If $V_\e$ is the $\e$-neighbourhood of the identity, then for every
$g\in V_\e$ and each $x\in \tilde X_n=X^{m_n}_n$ 
one has $d_1(g\cdot x, x)<\e\cdot a_n$, where $a_n=\diam X_n$. 
Consequently, if $g\in V_{\e/a_n}$, then $d_1(g\cdot x, x)<\e$.

Now let us turn to the group topology induced from $\Iso(\Ur)$.
Let 
\[V[x_1,\ldots,x_t;\e]=\{f\in\Iso(\Ur)\colon \forall i=1,2,\ldots,n,~~
d_{\Ur}(x_i,f(x_i))<\e\}
\]
be a standard neighbourhood of the identity
in $\Iso(\Ur)$. Here one can assume without loss in generality
that $x_i\in \cup_{n=1}^{\infty} X_n$, $i=1,2,\ldots,t$, because
the union of $X_n$'s is everywhere dense in $\Ur$. 
Let $k\in\N$ be such that $x_1,x_2,\ldots,x_t\in X_k$.
For all $n\geq k$, if $A\subseteq \tilde G_n$
contains at least half of all elements, the set 
$V_{\e/a_n}A$ is of Haar measure (taken in $\tilde G_n$)
at least $1-2 e^{-m_n \e/8 a_n^2}$,
according to Theorem \ref{th:product}. The set 
$V[x_1,\ldots,x_t;\e]\cdot A$ contains $V_{\e/a_n}A$ and so the measure of its
intersection with $\tilde G_n$ is at least as big.
According to the choice of numbers $m_n$ (Eq. \ref{eq:growth}), 
\[\mu_n(\tilde G_n\cap (V[x_1,\ldots,x_t;\e]\cdot A))\geq 1-e^{-n\e^2}. \]
By Proposition \ref{p:va}, the family of groups $\tilde G_n$
is L\'evy.
\qed\end{pf}

\section*{Acknowledgements}

I am grateful to the Organizers of the 6-th Iberoamerican Conference on
Topology and its Applications, held in charming Puebla from July 4--7, 2005,
for giving me an opportunity to deliver an invited lecture 
on the Urysohn space. This article is based on selected fragments of my lecture.
The research was supported by the Univesity of Ottawa internal research
grants and by the 2003-07 NSERC discovery grant 
``High-dimensional geometry and topological transformation groups.''


\begin{thebibliography}{100}

\bibitem{Cam} P. Cameron, \textit{The random graph,} in:
The Mathematics of Paul Erdos (J. Ne\v setril, R. L. Graham, eds.),
Springer, 1996, pp. 331-351.

\bibitem{Fre} M. Fr\'echet, \textit{Les espaces abstraits,} Paris, 1928.

\bibitem{GP2} 
T. Giordano and V. Pestov, \textit{Some extremely amenable groups related to
operator algebras and ergodic theory,} ArXiv e-print math.OA/0405288.

\bibitem{Gl1} S. Glasner,
\textit{On minimal actions of Polish groups,} Top. Appl. \textbf{85}
(1998), 119--125.

\bibitem{GTW} E. Glasner, B. Tsirelson, and B. Weiss,
\textit{The automorphism group of the Gaussian measure cannot act pointwise,}
ArXiv e-print math.DS/0311450, to appear in Israel J. Math.

\bibitem{GraL}
E. Granirer and A.T. Lau, \textit{Invariant means on locally compact groups,}
Ill. J. Math. \textbf{15} (1971), 249--257.

\bibitem{Gr} M. Gromov, \textit{Metric Structures for Riemannian and
Non-Riemannian Spaces,} Progress in Mathematics \textbf{152}, Birkhauser
Verlag, 1999.

\bibitem{GrM}
M. Gromov and V.D. Milman,
\textit{A topological application of the isoperimetric inequality,}
Amer. J. Math. \textbf{105} (1983), 843--854.

\bibitem{Gru} K.W. Gruenberg, \textit{Residual properties of infinite 
soluble groups,} Proc. London Math. Soc. (3), \textbf{7} (1957), 29--62.

\bibitem{HM} S. Hartman and J. Mycielski,  
{\it On the imbedding of topological groups into 
connected topological groups,}
Colloq. Math. {\bf 5} (1958), 167--169.

\bibitem{HC}  W. Herer and J.P.R. Christensen,
\textit{On the existence of pathological
submeasures and the construction of exotic topological groups,}
Math. Ann. \textbf{213} (1975), 203--210.

\bibitem{hui}  G.E. Huhunai\v svili,  
\textit{On a property of Uryson's universal metric space,}
Dokl. Akad. Nauk SSSR (N.S.)  \textbf{101} (1955), 607--610
(in Russian).

\bibitem{Kat} 
M. Kat\v{e}tov, {\it On universal metric spaces,} Gen.
Topology and its Relations to Modern Analysis and Algebra VI, Proc.
Sixth Prague Topol. Symp. 1986, Z. Frol\'ik (Ed.), Heldermann
Verlag, 1988, 323-330.

\bibitem{KPT} A.S. Kechris, V.G. Pestov and S. Todorcevic,
\textit{Fra\"\i ss\'e limits, Ramsey theory, and topological dynamics of
automorphism groups,} Geom. Funct. Anal. \textbf{15} (2005), 
106--189.

\bibitem{L}
M. Ledoux, \textit{The concentration of measure phenomenon}. Math. Surveys and 
Monographs, \textbf{89}, Amer. Math. Soc., 2001.

\bibitem{MKS}  W. Magnus, A. Karrass, and D. Solitar, 
\textit{Combinatorial Group Theory,} reprint of the 1976 second edition,
Dover Publications, Mineola, NY, 2004.

\bibitem{Meg2}
M.G. Megrelishvili, 
\textit{Operator topologies and reflexive representability,} in:
Nuclear groups and Lie groups (Madrid, 1999), 
Res. Exp. Math., \textbf{24}, Heldermann, Lemgo, 2001, pp. 197--208.

\bibitem{Meg} M.G. Megrelishvili,
\textit{Every semitopological semigroup compactification of the group
$H\sb +[0,1]$ is trivial}, Semigroup Forum \textbf{63} (2001), 357--370.

\bibitem{vM} J. van Mill, 
\textit{The Infinite-Dimensional Topology of Function Spaces,}
North-Holland Mathematical Library, \textbf{64},
North-Holland, Amsterdam, 2001. 

\bibitem{M86}  V.D. Milman, \textit{The concentration phenomenon and linear
structure of finite-dimensional normed spaces,} Proc. Intern. Congr. Math.
Berkeley, Calif., USA, 1986, pp. 961--975.

\bibitem{M00} V. Milman, \textit{Topics in asymptotic geometric analysis,}
Geom. Funct. Anal. --- special volume GAFA2000, 792--815.

\bibitem{MS}
V.D. Milman and G. Schechtman, \textit{Asymptotic theory of finite-dimensional 
normed spaces (with an Appendix by M. Gromov)}, 
Lecture Notes in Math., \textbf{1200}, Springer, 1986.

\bibitem{P99} 
V. Pestov, \textit{Topological groups:  where to from here?,} Topology
Proceedings {\bf 24} (1999), 421--502.

\bibitem{P02} V. Pestov, \textit{Ramsey--Milman phenomenon, 
Urysohn metric spaces,
and extremely amenable groups,} Israel Journal of Mathematics
\textbf{127} (2002), 317-358. \textit{Corrigendum,}  
ibid., \textbf{145} (2005), 375-379.

\bibitem{P05}
V. Pestov, \textit{Dynamics of infinite-dimensional groups and 
Ramsey-type phenomena,}
Publica\c c\~oes dos Col\'oquios de Matem\'atica, 
IMPA, Rio de Janeiro, 2005.

\bibitem{P05b}
V. Pestov, \textit{Forty annotated questions about large 
topological groups,} ArXiv e-print math.GN/0512564, December 2005,
11 pp., to appear in the second edition of Open 
Problems In Topology (Elliott Pearl, editor), Elsevier Science, 2007.

\bibitem{Rado} R. Rado, \textit{Universal graphs and universal functions,}
Acta Arithm. \textbf{9} (1954), 331--340.

\bibitem{Shtern} A. Shtern,
{\it Compact semitopological semigroups and reflexive representability
of topological groups,} Russian J. Math. Phys. {\bf 2} (1994),
131--132.

\bibitem{U25}
P. S. Urysohn, \textit{Sur un espace m\'etrique universel,}
C. R. Acad. Sci. Paris \textbf{180} (1925), 803--806.

\bibitem{U27} 
P. Urysohn, \textit{Sur un espace m\'etrique universel,} Bull. Sci.
Math. {\bf 51} (1927), 43--64 et 74--90.

\bibitem{U72} P.S. Urysohn, \textit{On the universal metric space,}
in: P.S. Urysohn. Selected Works (P. S. Alexandrov, ed.), vol. 2. 
Nauka, Moscow, 1972, pp. 747--769 (in Russian).

\bibitem{Usp86} V.V. Uspenskij,
{\it A universal topological group with countable base,}
Funct. Anal. Appl. {\bf 20} (1986), 160--161.

\bibitem{Usp90}
V.V. Uspenski\u\i,  {\it On the group of isometries of the
Urysohn universal metric space,} Comment. Math. Univ. Carolinae
{\bf 31} (1990), 181-182.

\bibitem{Usp98} V.V. Uspenskij, {\it On subgroups of minimal
topological groups,} 1998 preprint, later prepublished at
arXiv:math.GN/0004119.

\bibitem{Usp04} V. V. Uspenskij, \textit{The Urysohn universal metric space 
is homeomorphic to a Hilbert space,} Topology Appl. \textbf{139} (2004), 
145-149.

\bibitem{Ver98} A.M. Vershik, \textit{The universal Urysohn space, 
Gromov's metric triples, and random metrics on the series of positive 
numbers,} Russian Math. Surveys \textbf{53} (1998), 921--928.

\bibitem{Ver02} A.M. Vershik, \textit{Random metric space is Urysohn space,}
Dokl. Ross. Akad. Nauk \textbf{387} (2002), no. 6, 1--4.

\bibitem{Ver04} A.M. Vershik, \textit{The universal and random metric
spaces,} Russian Math. Surveys \textbf{356} (2004), 65--104.

\bibitem{Ver05} A.M. Vershik, \textit{Kantorovich metric:
initial history and little-known applications,} arXiv math.FA/0503035,
13 pp.

\bibitem{Ver05b} A.M. Vershik, \textit{Extensions of the partial isometries
of the metric spaces and finite approximation of the group of isometries of
Urysohn space,} preprint, 2005.

\end{thebibliography}
\end{document}